\documentclass[11pt]{article}

\usepackage{amsfonts,amssymb,graphicx}
\usepackage{amsfonts}
\usepackage{amsmath, amscd}
\usepackage{amstext}
\usepackage{amsthm}
\usepackage{xspace}

\newcommand{\Z}{\mathbb Z}
\newcommand{\N}{\mathbb N}
\newcommand{\C}{\mathbb C}

\newcommand{\R}{\mathbb R}
\newcommand{\Q}{\mathbb Q}

\newcommand{\virg}[1]{``#1"}
\newcommand{\vsni}{\vskip 0.2cm}

\newcommand{\Ps}{\rm {I\!\!\!P}}
\newcommand{\Es}{\rm {I\!\!\!E}}

\newtheorem {theorem}{Theorem}[section]

\newtheorem {lemma}[theorem]{Lemma}
\newtheorem {proposition}[theorem]{Proposition}
\newtheorem {corollary}[theorem]{Corollary}

\newtheorem {remark}[theorem]{Remark}

\theoremstyle{definition}

\numberwithin{equation}{section}

\def\be{\begin{equation}}
\def\ee{\end{equation}}
\newcommand{\Leq}[1]{\label{#1}\end{equation}}
\newcommand{\beqn}{\begin{eqnarray}}
\newcommand{\eeqn}{\end{eqnarray}}
\newcommand{\beqno}{\begin{eqnarray*}}
\newcommand{\eeqno}{\end{eqnarray*}}

\begin{document}

\title{Orderings of the rationals and dynamical systems}

\author{Claudio Bonanno
\thanks{Dipartimento di Matematica Applicata,
Universit\`a di Pisa, via F. Buonarroti 1/c, I-56127 Pisa, Italy,
email: $<$bonanno@mail.dm.unipi.it$>$} \and Stefano Isola
\thanks{Dipartimento di Matematica e Informatica, Universit\`a
di Camerino, via Madonna delle Carceri, I-62032 Camerino, Italy.
e-mail: $<$stefano.isola@unicam.it$>$}}

\maketitle

\begin{abstract}
\noindent
This paper is devoted to a systematic study of a class of binary trees encoding the structure of rational numbers both from arithmetic and dynamical point of view. The paper is divided into two parts. The first one is a critical review of rather standard topics such as Stern-Brocot and Farey trees and their connections with continued fraction expansion and the question mark function. In the second part we introduce a class of one-dimensional maps which can be used to generate the binary trees in different ways and study their ergodic properties. This also leads us to study some random processes (Markov chains and martingales) 
arising in a natural way in this context.

\end{abstract}
\vsni\vsni

\vsni
\vsni
\maketitle
\section{Part one: arithmetics}

\vsni
\noindent
{\bf Notational warning}: In the sequel we shall use the following notations: 
\begin{eqnarray}
I&:=& [0,1] \nonumber \\ J&:=&[0,\infty)\cup \{\infty\} \nonumber \\ 
\Q_1&:=&\Q \cap [0,1] \nonumber \\
\Q_p&:=&\left\{ \frac k {p^s} \, : \, s\in \N, \, 0\leq k \leq  p^s\right\} , \quad p\geq 2\nonumber
\end{eqnarray}

\subsection{A class of binary trees} 

We start with the Stern-Brocot (SB) tree $\cal T$, which is a way to order (and thus to count) the elements of $\Q^+$, the set of positive rational numbers, so that every number appears (and thus is counted) exactly once (see \cite{St}, \cite{Br} and, for a modern account, \cite{GKP}). 
The basic operation needed to construct  $\cal T$ is the
{\sl Farey sum}: given $\frac p q$ and $\frac {p'} {q'}$ in $\Q^+$ set  $${p\over q}\oplus {p'\over q'} = {p+p' \over q+q'}$$
One notes that the {\sl child} ${p\over q}\oplus {p'\over q'}$ turns out to be  in lowest terms whenever the  {\sl parents}
 $\frac p q$ and $\frac {p'} {q'}$ do. Moreover, the child always lies somewhere in between its parents, e.g., assuming  $\frac p q<\frac {p'} {q'}$, we have
 ${p\over q}\;< \; {p+p' \over q+q'} \;< \; {p'\over q'}$.

\vsni
\noindent
Starting from the {\sl ancestors} $0$ and $\infty$ 
(written `in lowest terms' ) one then writes genealogically one generation after the other using the above operation: 

\vskip 0.5cm

${0 \over 1}$ \hskip 10cm ${1\over 0}$

\vskip 0.2cm 
\hskip 5.17cm ${1\over 1}$

\vskip 0.5cm 
\hskip 2.5cm ${1\over 2}$ \hskip 4.95cm ${2\over 1}$

\vskip 0.3cm 
\hskip 1.15cm ${1\over 3}$ \hskip 2.3cm ${2\over 3}$  \hskip 2.3cm ${3\over 2}$  \hskip 2.4cm ${3\over 1}$

\vskip 0.3cm 
\hskip 0.5cm ${1\over 4}$ \hskip 0.95cm ${2\over 5}$  \hskip 1cm ${3\over 5}$  \hskip 1cm ${3\over 4}$
\hskip 0.9cm ${4\over 3}$ \hskip 1cm ${5\over 3}$ \hskip 1.05cm ${5\over 2}$  \hskip 1.cm ${4\over 1}$

\vskip 0.3cm 
\hskip 0.2cm ${1\over 5}$ \hskip 0.25cm ${2\over 7}$  \hskip 0.3cm ${3\over 8}$  \hskip 0.3cm ${3\over 7}$
\hskip 0.3cm ${4\over 7}$ \hskip 0.32cm ${5\over 8}$  \hskip 0.31cm ${5\over 7}$  \hskip 0.3cm ${4\over 5}$
\hskip 0.3cm ${5\over 4}$ \hskip 0.33cm ${7\over 5}$  \hskip 0.32cm ${8\over 5}$  \hskip 0.32cm ${7\over 4}$
\hskip 0.32cm ${7\over 3}$ \hskip 0.33cm ${8\over 3}$  \hskip 0.33cm ${7\over 2}$  \hskip 0.33cm ${5\over 1}$

\vskip 0.5cm 

\noindent
and so on. The easily verified property which makes the above interesting and useful is the following fact: if 
$\frac p q$ and $\frac {p'} {q'}$ are consecutive fractions at any stage of the construction then the unimodular relation $qp'-pq'=1$ is in force.

\vsni
\noindent
Finally, the subtree $\cal F$ of $\cal T$ having $\frac 1 2$ as root node and vertex set $\Q_1$ is called {\sl Farey tree}. It can be obtained exactly in the same way as $\cal T$ taking as ancestors $\frac 0 1$ and $\frac 1 1$ instead of $\frac 0 1$ and $\frac 1 0$. 

\begin{lemma}\label{moppa}
Let $\phi : J \to I$ be the
invertible map defined by $\phi (\infty)=1$ and
$$
\phi (x) = \frac x {x+1}\quad , \quad x\in \R^+
$$
Then 
$$
\phi ({\cal T}) = \cal F
$$
\end{lemma}
\noindent
{\sl Proof} It suffice to notice that $\phi(\frac 0 1 )=\frac 0 1$,  $\phi(\frac 1 0 )=\frac 1 1$, $\phi (\frac 1 1 ) = \frac 1 2$ and for $x, x' \in \Q^+$ we have $\phi(x)\oplus \phi(x')=\phi(x\oplus x')$. $\qed$

\vsni
\noindent
Another structure we shall deal with is the {\sl dyadic} tree $\cal D$, whose first two levels are as in $\cal F$ and then can be constructed from the root node  $\frac 1 2$ by  writing under each vertex $\frac p q$ the pair
$\frac {2p-1} {2q}$ and $\frac {2p+1} {2q}$. The vertex set of $\cal D$ is $\Q_2$. We shall see later how it is related to $\cal T$ and $\cal F$.

\subsection{Continued fractions and the $\{L,R\}$ coding} 

\noindent
Every  $x\in \Q^+$ 
appears exactly once in the above construction and corresponds to a
unique finite path on $\cal T$ starting at the root node $1\over 1$ and whose number of vertices equals the {\sl depth} of $x$, i.e. the level of $\cal T$ it belongs to. For $x\in \Q_1$ one may just consider  the path on the subtree
$\cal F$ which starts at the root node $1\over 2$ and whose number of vertices is the {\sl rank} of $x$. For $x\in \Q^+$ we have
\be
{\rm depth}(x) = [x] + {\rm rank}(\{x\})+1
\ee
In order to properly code these paths we start recalling
that every rational number $x\in \Q^+$ has a unique finite
continued fraction expansion \cite{Kh} 
$$ 
x = a_0+{1\over \displaystyle a_1 + 
{1\over {\ddots  \;  { \atop {\displaystyle + \frac{1}{ a_n} }}}}}\equiv [a_0; a_1,  \dots, a_n]
$$
with $a_0\geq 0$, $a_i\geq 1$ for $1\leq i <n$ and $a_n>1$.

\begin{lemma} Let $x\in \Q^+$ then
$$
\qquad x=[a_0;a_1,\dots ,a_n] \quad \Longrightarrow \quad {\rm depth}(x)= \sum_{i=0}^n a_i
$$
\end{lemma}
\noindent
{\sl Proof.} Setting
${\rm depth}({0\over 1})={\rm depth}({1\over 0})=0$ we have ${\rm depth}({1\over 1})=1$. 
Let now $x=[a_0;a_1,\dots ,a_n] $ be such that ${\rm depth}(x)=d>1$. Then, in order to reach the leaf $x$ from the root ${1\over 1}$
one has first to move $a_0$ steps to the right, thus reaching the node $a_0+{1\over 1}$. Then, moving $a_1$ steps to the left
one reaches $a_0+{1\over a_1+{1\over 1}}. \, $ $a_2$ further steps to the right reach the point $a_0+{1\over a_1+{1\over a_2+{1\over 1}}}$ and so on. In this way, one sees that the path to reach $x$ makes
exactly $n$ turns and the length of the blocks within the $(i-1)$-st and the $i$-th turn is given by the partial quotient $a_i$ for $1\leq i<n$, whereas the last block has length $a_n-1$.
More precisely, the blocks moving to the left are related to partial quotients with odd index, those moving to the right to those with even index. It then follows at once that $d=\sum_{i=0}^n a_i$. $\qed$

\vsni
\noindent
The argument sketched above actually allows us to say more. To this end, we shall first construct a 
matrix representation of the positive rationals.
We start noting that a given
$x \in \Q^+$ can be uniquely decomposed as
$$
x
={p\over q}\oplus {p'\over q'}\quad \hbox{with}\quad qp'-pq'=1  
$$  
The neighbours ${p\over q}$ and ${p'\over q'}$ are thus the parents
of  $x$ as an element of  $\cal T$. We then identify
$$
x \Longleftrightarrow 
\left(  \begin{array}{cc}
  p'&p\\
 q'&q \\
  \end{array}\right) \in SL(2,\Z)
$$
Note that the left column bears on the right parent and
viceversa. In this way, the root
node yields the identity matrix: $$ {1\over 1}={0\over 1}\oplus {1\over 0}
\Longleftrightarrow 
\left(  \begin{array}{cc}
  1&0\\
  0&1 \\
  \end{array}\right) $$
Moreover, given $M\in
SL(2,\Z)$ which represents the fraction $x \in \Q^+$, the matrix
$U\,M\, U$ represents the symmetric fraction $1/x$, with
$$
U=U^{-1}=
\left(  \begin{array}{cc}
  0&1\\
  1&0 \\
  \end{array}\right)
$$
In particular 
$$ {1\over 2} \Longleftrightarrow 
\left(  \begin{array}{cc}
  1&0\\
  1&1 \\
  \end{array}\right)=:L \quad \hbox{and} \quad
   {2\over 1} \Longleftrightarrow 
\left(  \begin{array}{cc}
  1&1\\
  0&1 \\
  \end{array}\right)=:R
$$ 
More generally, for $k\in \N$, 
$$ L^k= 
\left(  \begin{array}{cc}
  1&0\\
  k&1 \\
  \end{array}\right) \Longleftrightarrow \frac{1}{k+1}\quad\hbox{and}\quad R^k=\left(  \begin{array}{cc}
  1&k\\
  0&1 \\
  \end{array}\right) \Longleftrightarrow k+1
$$ 
Now, the point  $x$ considered above has in turn a
unique pair of (left and right) children, given by 
$$ {m\over
s}\oplus {m+n\over s+t} \quad\hbox{and}\quad {m+n\over s+t}\oplus
{n\over t} 
$$ respectively. 
Moreover,
$$
\left(  \begin{array}{cc}
  n&m\\
  t&s \\
  \end{array}\right) \left(  \begin{array}{cc}
  1&0\\
  1&1 \\
  \end{array}\right) 
  =\left(  \begin{array}{cc}
  m+n&m\\
 s+t&s \\
  \end{array}\right) 
 \Longleftrightarrow {m\over s}\oplus {m+n\over s+t}
$$ and 
$$ 
\left(  \begin{array}{cc}
  n&m\\
  t&s \\
  \end{array}\right)
  \left(  \begin{array}{cc}
  1&1\\
  0&1 \\
  \end{array}\right)=
  \left(  \begin{array}{cc}
  n&m+n\\
  t&s+t \\
  \end{array}\right)
\Longleftrightarrow {m+n\over
s+t}\oplus {n\over t} 
$$ 
In other words, the matrices $L$ and $R$,
when acting from the right, move to the left  and right child in
$\cal T$, respectively. Together with the argument of the proof given above this
yields the following result.
\begin{proposition}\label{codaggio}
To each entry $x\in {\cal T}$ there corresponds a unique element
$X\in SL(2,\Z)$, for which we have the following two possibilities:
\begin{itemize}
\item $x=[a_0;a_1,\dots ,a_n] $, $n$ even   $\Rightarrow$ $X=R^{a_0}L^{a_1}\cdots L^{a_{n-1}}R^{a_n-1}$
\item $x=[a_0;a_1,\dots ,a_n] $, $n$ odd $\Rightarrow$ $X=R^{a_0}L^{a_1}R^{a_2}\cdots L^{a_n-1}$
\end{itemize} 
\end{proposition}
\noindent
As an easy consequence we have the
\noindent
\begin{corollary}\label{path} Let $x=[a_0; a_1, \dots, a_n]$ with $a_n>1$ and $n$ even. Then its left and right children
in $\cal T$ are given by $x'=[a_0;a_1, \dots, a_n-1,2]$ and $x''=[a_0; a_1,
\dots, a_n+1]$, respectively. If instead $n$ is odd the expansions
for $x'$ and $x''$ are interchanged.
\end{corollary}
\noindent
{\sl Proof.} For $n$ even and larger than one we have
$
x=[a_0; a_1, \dots, a_n]
\Longleftrightarrow X=R^{a_0}L^{a_1}R^{a_2}\cdots R^{a_n-1}$.
Therefore
$
x' \Longleftrightarrow X'=R^{a_0}L^{a_1}R^{a_2}\cdots R^{a_n-1}L$ and $x''\Longleftrightarrow X''=R^{a_0}L^{a_1}R^{a_2}\cdots R^{a_n}$
which yield the claim. A similar reasoning applies for $n=0$ and $n$ odd.
$\qed$

\subsection{The infinite coding}\label{inf}
One can extend the above construction by associating to each $x\in \R^+$ a unique {\sl infinite} path 
in $\cal T$,
or else a unique semi-infinite word in $\pi (x) \in \{ L,R \}^\N$, in the natural way.  First, to
 $x\in \R^+\setminus \Q^+$ with infinite continued fraction
expansion $x=[a_0;a_1,a_2,a_3,\dots]$ there will correspond the (unique) sequence
$\pi (x) = R^{a_0}L^{a_1}R^{a_2}L^{a_3}\dots$, where now $R$ and $L$ are nothing but elements of a binary alphabet. 
For rational $x$ we can proceed as follows.
First we set $\pi (\frac 0 1 ) =L^\infty$ and 
$\pi(\frac 1 0 )=R^\infty$. Then note that
each $x\in \Q^+$ has two infinite paths which agree down to node
$x$: they are those starting with the finite sequence coding the
path to reach $x$ from the root node according to Proposition \ref{codaggio} and terminating with either
$RL^\infty$ or $LR^\infty$. We shall agree that $\pi (x)$
terminates with $RL^\infty$ or $LR^\infty$ according whether the
number of its partial quotients of $x$ is even or odd. Summarizing we have the following coding
\begin{itemize}
\item $x=[a_0;a_1,\dots ,a_n] $, $n$ even  $\Rightarrow$ $\pi(x)=R^{a_0}L^{a_1}\cdots R^{a_n}L^\infty$
\item $x=[a_0;a_1,\dots ,a_n] $, $n$ odd $\Rightarrow$ $\pi(x)=R^{a_0}L^{a_1}\cdots L^{a_n}R^\infty$
\item $x=[a_0;a_1,a_2,a_3,\dots]$ $\Rightarrow$
$\pi (x) = R^{a_0}L^{a_1}R^{a_2}L^{a_3}\dots$.
\end{itemize} 
\noindent
One easily checks that if $\succ$ denotes the lexicographic order on $\{L,R\}^\N$
then
$$
x>y \Longrightarrow \pi(x) \succ \pi(y)
$$
Finally, from the above it follows that for an  irrational $x$ the infinite path on $\cal T$ converging to $x$ coincides with the {\it slow continued fraction algorithm} (see, e.g., \cite{AO}).

\subsection{The (extended) question mark function}
\vskip 0.5cm
Given a
number $x\in \R^+$ with continued fraction expansion
$x=[a_0; a_1,a_2,\dots]$, one may ask what is the number obtained by
interpreting the sequence $\pi (x)$ defined in Section \ref{inf} as the
binary expansion of a real number in $(0,1)$. The number so obtained,
denoted $\rho (x)$, writes 
\be\label{defi}
\rho (x)=0\; .\;{\underbrace
{11\dots 1}_{a_0}}\, {\underbrace {00\dots 0}_{a_1}}\;{\underbrace
{11\dots 1}_{a_2}}\; \cdots \ee
 or, which  is the same,
\be
\rho (x)=1-\sum_{k\geq 0}(-1)^{k}\, 2^{-(a_0+\cdots
+a_k)}
\ee 
For instance $\rho (1/n)=1/2^{n}$ and $\rho (n)=1-1/2^n$ for all $n\geq 1$.
Setting $\rho (0)=0$ and 
$\rho (\infty )=1$ we see that $\rho : \R^+ \to I$  satisfies
\be\label{?}
\rho (x) = ?\circ \phi \,(x) 
\ee
where $\phi : J\to I$ is the map defined in Lemma \ref{moppa} and $?:I\to I$ is the {\sl Minkowski question mark function} \cite{M}, which
for $x=[0;a_1, a_2, \dots]$ is given by
\be\label{???}
? (x)=0\; .\;{\underbrace
{00\dots 0}_{a_1-1}}\, {\underbrace {11\dots 1}_{a_2}}\;{\underbrace
{00\dots 0}_{a_3}}\; \cdots 
\ee
Differently said, for $x\in (0,1)$ the number $?(x)$ is obtained by
interpreting the symbolic sequence corresponding to the path which starts from the root node $\frac 1 2$ and approaches $x$ along the Farey tree $\cal F$ as a 
binary expansion of a real number in $(0,1)$.

\vsni
\noindent
We now need a simple lemma.
\begin{lemma}\label{reci}
$$ 
x=[0;a_1,a_2, \dots ] \;\Longleftrightarrow \; 1-x=\left\{
\begin{array}{cl} 
[0;1+a_2,a_3, \dots ] & {\rm if}\ a_1=1  \\[0.5cm]
[0;1,a_1-1,a_2, \dots ] & {\rm if}\ a_1>1
\end{array} \right.
$$
\end{lemma}
\noindent
{\sl Proof.} Upon application of the identity
$$
\frac{1}{a+\frac 1 b}+\frac{1}{1+\frac {1} { a-1 + \frac 1 b}}=1\quad \hbox{and}\quad 
\frac{1}{1+\frac 1 {b+c}}+\frac{1}{1+b+c}=1\qquad \qed
$$

\begin{proposition}
The functions  $?$ and $\rho$ satisfy the functional equations 
\be
?(x)+?(1-x)=1\quad , \quad x\in I
\ee
and
\be
\rho (x) +\rho(1/x)=1\quad , \quad x\in J
\ee
respectively.
\end{proposition}
\noindent
{\sl Proof.} The equation for $?$ follows at once from Lemma \ref{reci} and (\ref{???}). That for $\rho$ 
then follows  by (\ref{?}). $\qed$

\vsni
\noindent
Additional properties of $\rho$ are inherited via (\ref{?}) from the following properties of $?$ (see \cite{Sa}, \cite{Ki}, \cite{VPB}, \cite{V}):

\begin{itemize}
\item $?(x)$ is strictly increasing from $0$ to $1$ and H\"older continuous of order $\beta ={\log
2\over \sqrt{5}+1}$
\item $x\in \Q_1$ iff $? (x)\in \Q_2$
\item $x$ is a quadratic irrational iff $? (x)$ is a (non-dyadic)
rational
\item $? (x)$ is a singular function: its derivative vanishes Lebesgue-almost everywhere
\end{itemize}

\begin{figure}[h]
\begin{center}
\includegraphics[width=10.0cm]{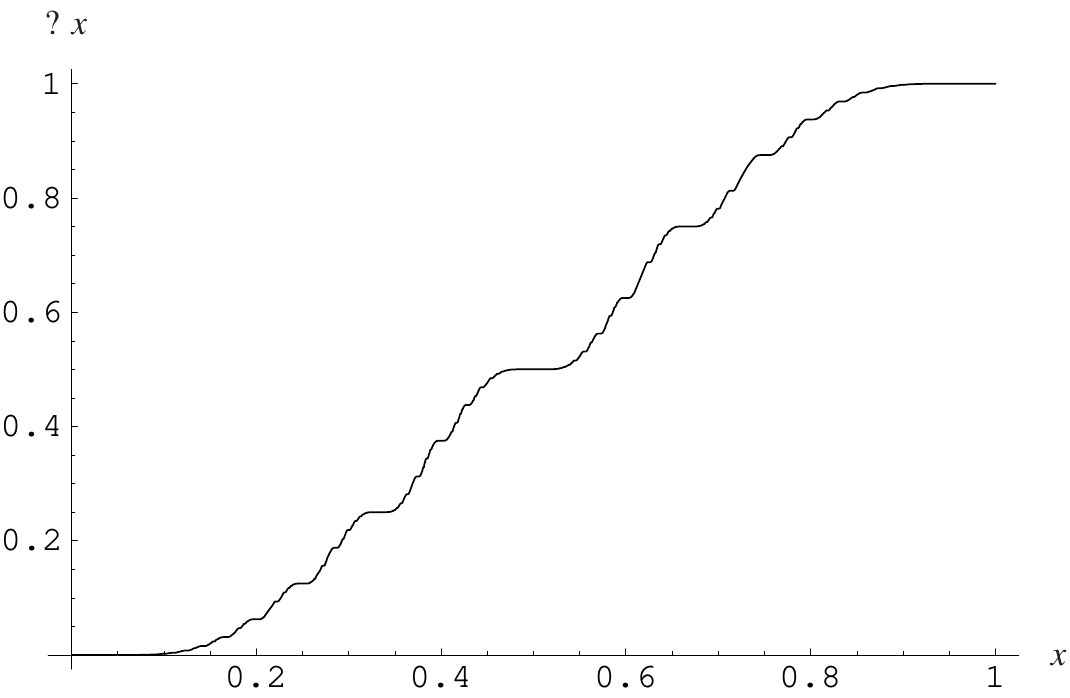}
\caption{The ? function}
\end{center}
\end{figure}

\vsni 
\noindent 
For any pair $\frac p q$ and $\frac {p'} {q'}$ of consecutive
fractions in $\cal T$ the function $\rho$ (as $?$ on $\cal F$) equates their child to the arithmetic
average. For instance we have
\be 
\rho\left({p+p'\over q+q'}\right) = {1\over 2} \left[
\rho \left({p\over q}\right) +\rho \left({p'\over q'}\right) \right] \ee 

\noindent
Therefore the functions  $\rho$ and $?$ map the SB tree $\cal T$ and the Farey tree $\cal F$, respectively, to
the {\sl dyadic tree} $\cal D$ mentioned above. Note that the set ${\cal D}_k$ of dyadic fractions belonging to the first $k+1$ levels of  $\cal D$ is the uniformly spaced sequence $l/2^{k}$, $l=0,1,\dots ,2^{k}$. Reducing to the lowest terms we get
$$
{\cal D}_0=\left({0\over 1}, {1\over 1}\right), \quad {\cal
D}_1=\left({0\over 1},  {1\over 2}, {1\over 1}\right),\quad {\cal
D}_2=\left({0\over 1},  {1\over 4}, {1\over 2},  {3\over 4}, {1\over
1}\right)
$$ and so on. Hence, an immediate consequence of
the fact that $\rho (\cal T)=?(\cal F) = {\cal D}$ is that $\rho (x)$ and $?(x)$ are the asymptotic
distribution functions of the sequences of SB fractions and Farey fractions, respectively. 
\begin{theorem}\label{distri}  Set ${\cal T}_{k}:=\{ \frac p q \in {\cal T}\, : \, {\rm depth}(x)\leq k\}$.
Since $$ \left| x - {\# \{ {p\over q}\in {\cal
D}_{k}\, : \, {p\over q}\leq x\} \over 2^k} \right| \leq 2^{-k} $$ we have $$
\left| \rho (x) - {\# \{ {p\over q}\in {\cal T}_{k}\, : \, {p\over q}\leq x\} \over 2^k} \right| \leq 2^{-k}
 $$
The same holds for $?$ with ${\cal T}_{k}$ replaced by ${\cal F}_{k}:=\{ \frac p q \in {\cal F}\, : \, {\rm rank}(x)\leq k\}$.
\end{theorem}
\noindent
In particular, the Fourier-Stieltjes
coefficients of $\rho$ and $?$ are as in the following
\begin{corollary} \label{fourier}
Let
$$ c_n:= \int_0^\infty e^{2\,\pi\, i\, n\, x} d\rho (x)
$$ then $$c_n
=\lim_{k \to \infty} {1\over 2^{k}}\sum_{{p\over q}\in {\cal
T}_{k}} e^{2\,\pi\, i\, n\, {p\over q}}
$$
The same holds for the coefficients of $?$ with ${\cal T}_{k}$ replaced by ${\cal F}_{k}$.
\end{corollary}
\noindent

\subsection{Permuted trees} 

\vskip 0.5cm
Let $X \in SL\, (2,\Z)$ represent a number $x\in \Q^+$ as above and write it as $X=\prod_{i=1}^k
M_i$ where
$M_i\in
\{L,R\}$ and $k={\rm depth}(x)$. Let us denote by $\hat x$ the positive rational number represented by the reversed matrix product 
${\hat X}=\prod_{i=k}^1 M_i$. Clearly ${\rm depth}(x)={\rm depth}(\hat x)$ but $x=\hat x$ if and only if the sequence $M_1\dots
M_k$ is a palindrome. The permutation map $x \mapsto \hat x$ yields the permuted version $\cal {\hat T}$ of the SB tree whose first five levels are shown below (the ancestors $\frac 0 0$ and $\frac 1 0$ are omitted). 

\vskip 0.5cm
 
\hskip 5.17cm ${1\over 1}$

\vskip 0.3cm 
\hskip 2.4cm ${1\over 2}$ \hskip 5cm ${2\over 1}$

\vskip 0.3cm 
\hskip 1.1cm ${1\over 3}$ \hskip 2.3cm ${3\over 2}$  \hskip 2.3cm ${2\over 3}$  \hskip 2.4cm ${3\over 1}$

\vskip 0.3cm 
\hskip 0.5cm ${1\over 4}$ \hskip 0.9cm ${4\over 3}$  \hskip 1cm ${3\over 5}$  \hskip 1cm ${5\over 2}$
\hskip 0.9cm ${2\over 5}$ \hskip 1cm ${5\over 3}$ \hskip 1.1cm ${3\over 4}$  \hskip 1cm ${4\over 1}$

\vskip 0.3cm 
\hskip 0.2cm ${1\over 5}$ \hskip 0.25cm ${5\over 4}$  \hskip 0.3cm ${4\over 7}$  \hskip 0.3cm ${7\over 3}$
\hskip 0.3cm ${3\over 8}$ \hskip 0.32cm ${8\over 5}$  \hskip 0.31cm ${5\over 7}$  \hskip 0.3cm ${7\over 2}$
\hskip 0.3cm ${2\over 7}$ \hskip 0.33cm ${7\over 5}$  \hskip 0.32cm ${5\over 8}$  \hskip 0.32cm ${8\over 3}$
\hskip 0.32cm ${3\over 7}$ \hskip 0.33cm ${7\over 4}$  \hskip 0.33cm ${4\over 5}$  \hskip 0.33cm ${5\over 1}$

\vskip 0.8cm

\begin{lemma}\label{desce}  Under the ancestors $\frac 0 1$ and $\frac 1 0$, the permuted SB tree $\hat {\cal T}$ can be constructed starting from the root node $\frac 1 1$ and 
writing under each vertex $\frac p q$ the set of descendants $\{\frac p {p+q}, \frac {p+q} {q}\}$. \end{lemma}
\noindent
{\sl Proof.} Note that 
$$
\left(  \begin{array}{cc}
  1&0\\
  1&1 \\
  \end{array}\right) \left(  \begin{array}{cc}
  n&m\\
  t&s \\
  \end{array}\right) 
  =\left(  \begin{array}{cc}
  n&m\\
 n+s&m+t \\
  \end{array}\right) 
 \Longleftrightarrow {m+n\over m+n+s+t}
 $$ and 
$$
\left(  \begin{array}{cc}
  1&1\\
  0&1 \\
  \end{array}\right) \left(  \begin{array}{cc}
  n&m\\
  t&s \\
  \end{array}\right) 
  =\left(  \begin{array}{cc}
   n+s&m+t\\
 s&t \\
  \end{array}\right) 
 \Longleftrightarrow {m+n+s+t\over s+t}
 $$
In other words, the matrices $L$ and $R$,
when acting from the left give but the left and right descendants, respectively. 

\noindent
Also note that if ${p\over q}=[a_0;a_1, \dots, a_n]$ then ${p\over
p+q}=[0;1,a_0,a_1, \dots, a_n]$ and ${p+q\over q}=[a_0+1; a_1, \dots, a_n]$.
Therefore $${\rm depth}({q\over p+q})={\rm depth}({p+q\over q})={\rm
depth}({p\over q})+1$$This yields the claim. $\qed$

\vsni
\noindent
\begin{remark}
The tree ${\hat {\cal T}}$ has been considered in \cite{CW} where the authors argued that if we read it row by row, and each row from left to right, then for $i\geq 2$ we can write the $i$-th element in the form $x_{i}=b(i-2)/b(i-1)$, where $b(n)$ is the number of {\rm hyperbinary} representation of $n$, namely the number of ways of writing the integer $n$ as a sum of powers of two, each power being used at most twice. For example $8=2^3=2^2+2^2=2^2+2+2=2^2+2+1+1$ and therefore $b(8)=4$. 
This property plainly entails that when reading from left to right any sequence of fractions with fixed depth the denominator of each fraction is the numerator of its successor. 
\end{remark}
\vsni
\noindent
We finally define the corresponding permutation of both the Farey tree $\cal F$ and the dyadic tree $\cal D$, denoted ${\hat {\cal F}}$ and ${\hat {\cal D}}$ 
 respectively. Clearly we have ${\hat {\cal F}}=\phi ({\hat {\cal T}})$ (see Lemma \ref{moppa}).  
 Reasoning as above one easily obtains the following simple genealogical rules:
 \begin{lemma}\label{genera}
  Under the ancestors $\frac 0 1$ and $\frac 1 1$, the permuted trees $\hat {\cal F}$ and ${\hat {\cal D}}$  can be constructed starting from the root node $\frac 1 2$ and 
writing under each vertex $\frac p q$ the sets of descendants $\{ \frac p {p+q}, \frac q {2q-p} \}$ and 
$\{ \frac p {2q}, \frac {p+q} {2q} \}$, respectively. \end{lemma}
\vsni
\noindent
The first five levels of ${\hat {\cal F}}$ are

\vskip 1cm 
\hskip 5.17cm ${1\over 2}$

\vskip 0.3cm 
\hskip 2.4cm ${1\over 3}$ \hskip 5.1cm ${2\over 3}$

\vskip 0.3cm 
\hskip 1.1cm ${1\over 4}$ \hskip 2.4cm ${3\over 5}$  \hskip 2.21cm ${2\over 5}$  \hskip 2.45cm ${3\over 4}$

\vskip 0.3cm 
\hskip 0.5cm ${1\over 5}$ \hskip 0.9cm ${4\over 7}$  \hskip 1cm ${3\over 8}$  \hskip 1.05cm ${5\over 7}$
\hskip 0.8cm ${2\over 7}$ \hskip 1.05cm ${5\over 8}$ \hskip 1.1cm ${3\over 7}$  \hskip 0.905cm ${4\over 5}$

\vskip 0.3cm 
\hskip 0.2cm ${1\over 6}$ \hskip 0.2cm ${5\over 9}$  \hskip 0.2cm ${4\over 11}$  \hskip 0.2cm ${7\over 10}$
\hskip 0.15cm ${3\over 11}$ \hskip 0.2cm ${8\over 13}$  \hskip 0.2cm ${5\over 12}$  \hskip 0.2cm ${7\over 9}$
\hskip 0.2cm ${2\over 9}$ \hskip 0.25cm ${7\over 12}$  \hskip 0.15cm ${5\over 13}$  \hskip 0.2cm ${8\over 11}$
\hskip 0.2cm ${3\over 10}$ \hskip 0.2cm ${7\over 11}$  \hskip 0.2cm ${4\over 9}$  \hskip 0.2cm ${5\over 6}$
\vskip 1cm
\noindent
and the corresponding levels of ${\hat {\cal D}}$ are

\vskip 1cm

\hskip 5.17cm ${1\over 2}$

\vskip 0.3cm 
\hskip 2.4cm ${1\over 4}$ \hskip 5.3cm ${3\over 4}$

\vskip 0.3cm 
\hskip 1.1cm ${1\over 8}$ \hskip 2.3cm ${5\over 8}$  \hskip 2.5cm ${3\over 8}$  \hskip 2.45cm ${7\over 8}$

\vskip 0.3cm 
\hskip 0.4cm ${1\over 16}$ \hskip 0.7cm ${9\over 16}$  \hskip 0.8cm ${5\over 16}$  \hskip 0.9cm ${13\over 16}$
\hskip 0.9cm ${3\over 16}$ \hskip 1cm ${11\over 16}$ \hskip 0.9cm ${7\over 16}$  \hskip 0.9cm ${15\over 16}$

\vskip 0.3cm 
\hskip 0.1cm ${1\over 32}$ \hskip 0.1cm ${17\over 32}$  \hskip 0.1cm ${9\over 32}$  \hskip 0.1cm ${25\over 32}$
\hskip 0.15cm ${5\over 32}$ \hskip 0.2cm ${21\over 32}$  \hskip 0.2cm ${13\over 32}$  \hskip 0.2cm ${29\over 32}$
\hskip 0.2cm ${3\over 32}$ \hskip 0.25cm ${19\over 32}$  \hskip 0.2cm ${11\over 32}$  \hskip 0.2cm ${27\over 32}$
\hskip 0.2cm ${7\over 32}$ \hskip 0.2cm ${23\over 32}$  \hskip 0.2cm ${15\over 32}$  \hskip 0.2cm ${31\over 32}$
\vskip 0.8cm 
\noindent

\vskip 0.5cm 

\subsection{Random walks on the permuted trees}

\noindent
We now construct a sequence of random variables 
$Z_1, Z_2, \dots $ on $\Q^+$  defined recursively in the following way: set $Z_1=\frac 1 1$ and if $Z_k= \frac p q$ then either $Z_{k+1}=\frac p {p+q}$ or $Z_{k+1}=\frac {p+q} q$, both with probability $\frac 1 2$. 
The sequence $(Z_k)_{k\geq 1}$ can be regarded as a (symmetric) random walk on ${\hat {\cal T}}$.

\begin{theorem} \label{rw} The random walk $(Z_k)_{k\geq 1}$ enters any non-empty open interval $(a,b)\subset \R^+$
almost surely.
\end{theorem}
\noindent
{\sl Proof.} Pick up  an irrational number $x=[a_0;a_1,\dots ]\in (a,b)$. Then for $n$ large enough we can find a closed subinterval $A \subset (a,b)$ such that the c.f. expansion of each element of $A$ starts as $[a_0;a_1, \dots , a_n, \dots]$. To fix the ideas and with no loss, let $n$ be odd. Then, according to the above and Proposition \ref{codaggio} the path on ${\hat {\cal T}}$ starting from the root $\frac 1 1$ and entering $A$ (for the first time) will eventually end with the word $W=L^{a_n-1} \cdots R^{a_2} L^{a_1} R^{a_0}$. Hence it has the form $UW$ with prefix $U\in \{L,R\}^*$ so that $UW$ does not contain subwords equal to $W$ but $W$ itself. We now proceed by induction on the length $\ell$ of $W$, calling $W_\ell$ the word of length $\ell$.
If $\ell =1$ then there is exactly one prefix $U$ of lenght $k$ for each $k\geq 1$ (e.g. if $W_1=L$ then $U=R^k$ is the only possible prefix)
occurring with probability $2^{-k}$. Summing over the prefixes we get $\sum_{k\geq 1} 2^{-k}=1$. Theferore the claim is true for $\ell =1$. Now suppose it is true for $\ell =m$. When passing to $\ell =m+1$ either $W_{m+1}=W_mL$ or $W_{m+1}=W_mR$, hence we have two families of paths
$UW_mL$ and $UW_mR$, one of which being $UW_{m +1}$ and thus, by the induction hypothesis, having probability $\frac 1 2$. We are now left with all paths starting with the `bad' ones and eventually ending with $W_{m+1}$. But then we can use the self-similarity of the tree and iterate the above construction. Suppose for instance that the `bad' set was  $UW_mR$, that is $W_{m+1}=W_mL$. Then at some point we will end up with the alternative  $UW_mRU'W_mL$ and $UW_mRU'W_mR$ for some $U'\in \{L,R\}^*$, and the `good' set $UW_mRU'W_mL$ has probability $\frac 1 2 \cdot \frac 1 2$. Iteration of this argument yields the
probability $ \frac 1 {2^2} +\frac 1 {2^3}+\cdots  = \frac 1 2$ which has to be added to the probability $\frac 1 2$ of the initial `good' set $UW_mL$. $\qed$

\begin{remark} The above result can be easily extended to both ${\hat {\cal F}}$ and ${\hat {\cal D}}$. However, it seems to be peculiar of the particular permutation which defines these trees, in particular it is plainly false for the
original Stern-Brocot tree ${{\cal T}}$ (as well as for $\cal F$ and $\cal D$). We shall see later a further generalisation.
\end{remark}
\section{Part two: dynamics}
\noindent
We shall now be dealing with a class transformations which generate the permuted trees ${\hat {\cal T}}$, ${\hat {\cal F}}$ and ${\hat {\cal D}}$, respectively, either one generation after the other or in genealogical way, i.e. producing elements with increasing depth.

\subsection{Rank one ergodic transformations with dense orbits of rationals}
\noindent
It was noticed in \cite{N} that the sequence $x_i$ of elements of $\hat {\cal T}$ satisfies the iteration
\be\label{itera}
x_{i+1}=\frac{1} { 1 -\{x_i\}+[x_i]}\quad ,\quad  i\geq 0
\ee
We are thus led to study the map\footnote{The study of this map was suggested  to one us us (C.B.) by Don Zagier}
$R: J \to  J$ given by $R(\infty):=0$ and
\be\label{map}
\qquad \qquad R(x):= \frac{1} { 1 -\{x\}+[x] }\quad ,\quad  x\in\R^+
\ee

 \begin{proposition}\label{uno}

\noindent
\begin{enumerate}

\item $R$ is an automorphism of $\R^+$
\item for any $x\in \R^+$, $R(x)\in \Q^+$ if and only if $x\in \Q^+$
\item $R$ counts the set $\Q^+\cup \{0\}\cup\{\infty\}$ in the following sense:
let $x_i$ be the sequence obtained  by reading  $\hat {\cal T}$ row by row and each row from left to right  (except the zero-th one), then
$x_i = R^i(\small{\frac{1}{0}})$.
\end{enumerate}
\end{proposition}

\noindent
{\sl Proof.}
One easily checks that $R$ is one-to-one and onto, with inverse
\be
R^{-1}(x) = 2n+1 -\frac 1 x \quad , \quad \frac{1}{n+1}\leq x < \frac 1 n \quad , \quad n\geq 0\, .
\ee
This proves 1. Statement 2 is immediate.  
Moreover, if $x\in \N$ then
$R(x)=1/(x+1)$ so that ${\rm depth}(R(x))={\rm depth}(x)+1$. If instead $x=[a_0;a_1,\dots ,a_n]$ with $n\geq 1$ then
$R(x)=1/(a_0+1-\{x\})$ so that ${\rm depth}(R(x))={\rm depth}(x)$ since ${\rm depth}(1/x)={\rm
depth}(x)$ and ${\rm rank}(\{x\})={\rm rank}(1-\{x\})$. This yields the first part of statement 3. To see the second part we start observing that if we write $x$ in the form $x=(kq+r)/q$ with $k\geq 0$ and $0\leq r <q$ we have $R(x)= q/(kq+q-r)$.  Now if $k=0$ then $x=r/q$ and $R(x)= q/(q-r)$, namely $x$ and $R(x)$ are left and  right descendants of the fraction $r(/q-r)$. If instead $k>0$ then $x$ is the right descendant of $x'=((k-1)q+r)/q$ whereas
$R(x)$ is the left descendant of  $x''=q/(kq-r)$, and $x''=R(x')$. $\qed$

\begin{remark}{\rm 
Note that, although the sequence $(x_i)_{i\geq 0}$ defined in (\ref{itera}) is dense in $\R^+$, it `diffuses' only logarithmically. Indeed we have
$x_i=n$ for $i=2^n$ and therefore $\sup_{0<i\leq n} x_i = {\cal O}(\log n)$. In fact, from what is proved below it follows that all 
orbits $\{ R^i(x) \, , \, i\geq 0\}$, $x\in \R^+$, are dense and have this property.
An automorphism of the unit circle with similar properties has been constructed in \cite{Bo}.}
\end{remark}
\noindent
We now restrict to the unit interval and consider two automophisms on it. The first one is the the map $S: I \to I$ defined by (see Lemma \ref{moppa})
\be\label{esse}
S(x):= \phi \circ R \circ \phi^{-1}(x)
\ee
or else by $S(1)=0$ and
\be\label{map}
S(x) = \frac{1} { 2 -\left\{  \displaystyle \frac x{1-x}\right\}+\left[\displaystyle \frac x{1-x}\right] }\quad , \quad  x\in [0,1)
\ee
Its inverse is
\be
S^{-1}(x) = \frac {2\, n\, x-1}{(2n+1) x-1} \quad , \quad \frac{1}{n+1}\leq x < \frac 1 n \quad , \quad n\geq 1
\ee
The second is the classical {\sl Von Neumann-Kakutani transformation} $T: I \to I$ given by $T(1):=0$ and 
\be
T(x):= x+ \frac 3 {2^{n+1}} -1 \quad , \quad 1-\frac 1 {2^n} \leq x < 1-\frac 1 {2^{n+1}} \quad , \quad n\geq 1
\ee
It was defined in \cite{VN} and is also called {\sl van der Corput's  transformation} or else {\sl dyadic rotation}. 

\begin{figure}[h]
\begin{center}
\includegraphics[width=13.0cm]{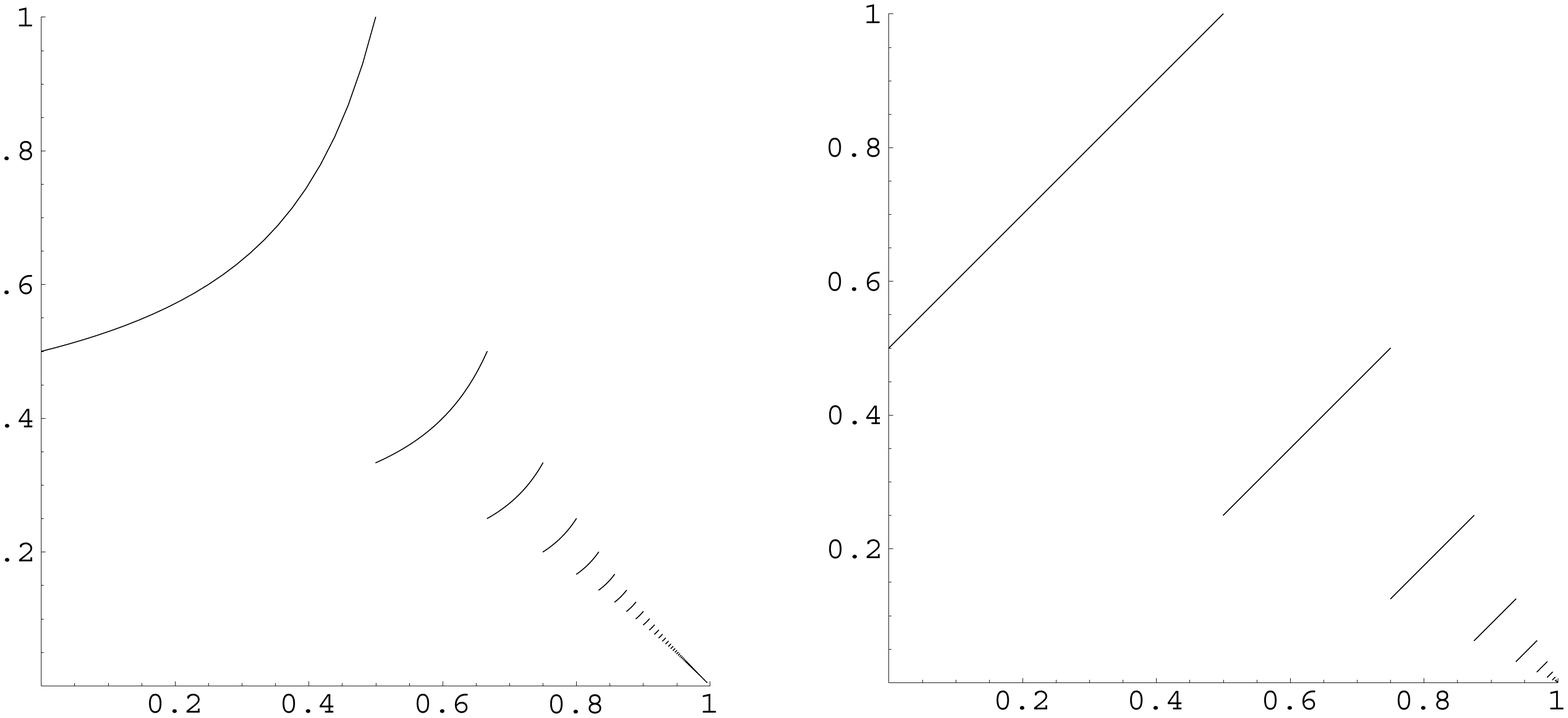}
\caption{{\small The maps $S$ (left) and  $T$ (right).}}
\label{indu}
\end{center}
\end{figure}

\begin{theorem} \label{conju} We have the following
commutative diagram
$$
\begin{CD}
J  @>\phi >>  I @> ? >> I \\
@VRVV  @VVSV @VVTV \\
J  @>\phi >>  I @> ? >> I \\
\end{CD}
$$
\end{theorem}
\noindent
{\sl Proof.} The left half follows immediately from (\ref{esse}). To see the right half part let $x\in (0,1)$
be given by $x=[0;a_1,a_2, \dots]$. We have $\frac x {1-x} = [0;a_1-1,a_2, \dots]$ so that
$$
\left[\frac x{1-x}\right] =\left\{
\begin{array}{cl}
a_2 & {\rm if}\ a_1=1 \\[0.5cm]
0  & {\rm if}\ a_1>1
\end{array} \right. \quad , \quad \left\{  \frac x{1-x}\right\}=\left\{
\begin{array}{cl}
[0;a_3,a_4, \dots] & {\rm if}\ a_1=1 \\[0.5cm]
[0;a_1-1,a_2, \dots]  & {\rm if}\ a_1>1
\end{array} \right. 
$$
and therefore 
\be \label{esse1}
S(x)=\left\{
\begin{array}{cl}
(2+a_2-[0;a_3,a_4, \dots] )^{-1}& {\rm if}\ a_1=1 \\[0.5cm]
(2-[0;a_1-1,a_2, \dots])^{-1}  & {\rm if}\ a_1>1
\end{array} \right. 
\ee
Using Lemma \ref{reci}  this becomes
\be \label{esse2}
S(x)=\left\{
\begin{array}{cl}
[0;a_2+1,1,a_3-1,a_4, \dots]  & {\rm if}\ a_1=1 \\[0.5cm]
[0;1,1,a_1-2,a_2, \dots]  & {\rm if}\ a_1>1
\end{array} \right. 
\ee
where for $a_k\geq \ell \geq 1$ for some $k\geq 1$ we set $$[a_0;a_1,a_2,\dots ,a_{k-1},a_k-\ell,a_{k+1}, \dots]=[a_0;a_1,a_2,\dots ,a_{k-1}+a_{k+1}, \dots]$$ if $a_k=\ell$.
On the other hand, the map $T(x)$ is also  named dyadic rotation because of the following fact (see, e.g. \cite{PF}, p.120): if we expand $x\in [0,1]$ in base two, i.e. we write $x=\sum_{k=0}^\infty \omega_k 2^{-k-1}=0.\omega_1\omega_2\dots$ with $\omega_k\in\{0,1\}$, it acts as  $T(0.111\dots)=0.000\dots$ and for $n\geq 1$
\be
T(0.\, {\underbrace {11\dots 1}_{n-1}}\,0\, \omega_{n+1}\omega_{n+2}\dots)=0.\, {\underbrace {00\dots 0}_{n-1}}\, 1\,\omega_{n+1}\omega_{n+2}\dots
\ee
Therefore if $x=[0;a_1,a_2, \dots ]$ then using (\ref{???}) we find
\be
T(?(x))=\left\{
\begin{array}{cl}
0.\,  {\underbrace {00\dots 0}_{a_2}}\, 1\,{\underbrace {00\dots 0}_{a_3-1}}\, {\underbrace {11\dots 1}_{a_4}}\dots
 & {\rm if}\ a_1=1 \\[0.5cm]
0.\, 1\, {\underbrace {00\dots 0}_{a_1-2}}\, \, {\underbrace {11\dots 1}_{a_2}}\dots & {\rm if}\ a_1>1
\end{array} \right. 
\ee 
which is identical to what we obtain applying $?$ to (\ref{esse2}). $\qed$

\vsni
\noindent 
We now derive some consequences from the above theorem.
\begin{corollary}\label{itera} The maps
$S$ and $T$ count the sets $\Q_1$ and $\Q_2$, respectively. More specifically, let for instance 
$y_i$ be the sequence obtained  by reading  $\hat {\cal F}$ row by row and each row from left to right  (except the zero-th one), then
$y_i = S^i(\small{\frac{1}{1}})$. A similar statement holds for $\hat {\cal D}$ and $T$.
\end{corollary}
\noindent
{\sl Proof.}  Follows at once from Proposition \ref{uno} and Theorem \ref{conju}. $\qed$

\begin{corollary}\label{uniergo} The systems $(J,R)$ and $(I,S)$ are uniquely ergodic and hence $(J,R, d\rho)$ and $(I,S,d?)$ are ergodic.
\end{corollary}
\noindent
{\sl Proof.} The first statement follows from the above as topological conjugacy preserves unique ergodicity and the system $(I,T)$ has this property. Moreover, the Lebesgue neasure $dx$  is $T$-invariant so that by the above and (\ref{?}) the maps $R$ and $S$ preserve the measures $d\rho$ and $d?$ respectively.$\qed$

\begin{corollary} The systems $(J,R)$ and $(I,S)$ are of rank one. Moreover they have the same spectrum which is discrete with eigenvalues
$e^{2\pi i \alpha}$ for any dyadic rational  $\alpha$.
\end{corollary}
\noindent
{\sl Proof.} 
The system $(I,T)$ has this property. Let us briefly recall how this is obtained. One start setting $A(1,n)=[0,2^{-n})$ for $n\geq 0$ and noticing that
$T$ maps in an affine way $A(i,n)=T^{i-1}A(1,n)$ onto $A(i+1,n)$ for $i=1,\dots , 2^n$. Clearly, these intervals are not ordered lexicographically but in the way induced by $T$. For example, for $n=3$ we have $000 \mapsto 100 \mapsto 010 \mapsto 110 \mapsto 001 \mapsto 101 \mapsto 011 \mapsto 111$. One may then write the so ordered intervals one above the other, thus making a stack which partitions the whole space.  The action of $T$ is then that of climbing up one level in the $n$-stack but is not defined on the top level. At step $n+1$, i.e. looking at the action of the iterates of $T$ on $A(1,n+1)$, the stack is cutted into two equal halves and the right half is stacked on the left half. This defines the action of $T$ on a finer partition of the space. This procedure eventually leads to the knowledge of $T$ on the whole space.
Finally, to get the same property for $(I,S)$ it will suffice to follow the above procedure with the family of intervals $B(i,n)=?^{-1}(A(i,n))$ (stacked in same order). Clearly, although all the intervals $A(i,n)$, $i=1, \dots , 2^n$ do have the same length $2^{-n}$, the corresponding $B(i,n)$ do not. A similar construction can be done for $(J,G)$ with the intervals $C(i,n)=\phi^{-1}(B(i,n))$.
The last assertion follows again from the same property for $(I,T)$ along with topological conjugacy (see, e.g., \cite{PF}, p.23). $\qed$

\vsni
\noindent 
Finally, setting $e_n(x):=e^{2\,\pi\, i\, n\, x}$, the Fourier-Stieltjes coefficients of $\rho$ are
\be
c_n:= \int_0^\infty e_n(x)  d\rho (x)
\ee
By Corollary \ref{fourier} they can be computed as
\be
 c_n
=\lim_{k \to \infty} {1\over 2^{k}}\sum_{{p\over q}\in {\cal
T}_{k}} e_n(\textstyle \frac p q)
\ee
On the other hand, the unique ergodicity of $(J,R)$ following from Corollary \ref{uniergo}
ensures that they can also be computed as ergodic means in a uniform way.
\begin{corollary} We have, uniformly  for $x \in J$,

$$
c_n=\lim_{N\to \infty}{1\over N} \sum_{k=0}^{N-1} e_n (R^k(x))$$
The same holds for the coefficients of $?$ with $R$ replaced by $S$.
\end{corollary}
\begin{remark}
An interesting question is whether $c_n\to 0$ as $n\to \infty$ (see \cite{Sa}). 
A solution to this question would give some insights into the degree of uniformity of the distribution of the dense sequence 
$\{R^k(x)\}_{k\geq 0}$, $x\in J$, and in particular of the permuted Stern-Brocot sequence obtained setting $x=1$.
\end{remark}
 \subsection{Markov maps and transfer operators}
\noindent
We now introduce three non-invertible maps which generate the trees $\hat {\cal T}$, $\hat {\cal F}$ and $\hat {\cal D}$ genealogically, i.e. via descendants. With the notations of Theorem \ref{conju}, the first one is the map $G: J \to J$ given by
\be\label{G} G(x)=\left\{
\begin{array}{cl}
\displaystyle {x\over 1-x}  & {\rm if}\ 0\leq x<1\\[0.5cm]
 x-1& {\rm if}\ x \geq 1
\end{array} \right.  
\ee
The second is the  {\sl modified Farey map} $F:I\to I$ given by
\be\label{farey} F(x)=\left\{
\begin{array}{cl}
\displaystyle {x\over 1-x}  & {\rm if}\ 0\leq x< {1\over 2} \\[0.5cm]
 2-{1\over x} & {\rm if}\ {1\over 2} \leq  x\leq 1
\end{array} \right.  
\ee
and the third is  {\sl doubling map} $D:I\to I$ given by
\be
D(x)=2x \, ({\rm mod}\, 1)
\ee
They are expansive  orientation preserving piecewise analytic endomorphisms such that the sets $G^{-1}(x)$, $F^{-1}(x)$ and $D^{-1}(x)$ are composed exactly by two points for each $x$. More specifically 
\begin{eqnarray} \label{inve}
G^{-1}(x)&=&\left\{ \frac x {1+x} , x+1 \right\}\quad , \quad  x\in J  \nonumber \\
F^{-1}(x)&=&\left\{ \frac x {1+x} , \frac 1 {2-x} \right\}  \quad , \quad
D^{-1}(x)=\left\{ \frac x 2 , \frac x  2 + \frac 1 2 \right\} \quad , \quad  x\in I\nonumber 
\end{eqnarray}
Both $F$ and $D$ fix the boundary points $0$ and $1$, but for  $F$ these are {\sl indifferent fixed points}, i.e.  $F'(0)=F'(1)=1$. More specifically, $0$ is a {\sl weakly repelling} fixed point whereas $1$ is {\sl weakly attracting}. 
On the other hand we can say that $G$ has two indifferent fixed points at $0$ and $\infty$. 
\begin{theorem} The permuted tree  $\hat {\cal T}$  can be constructed genealogically from its root  $\frac 1 1$ by
writing under each  leaf $x$ the set of descendants $G^{-1}(x)$. The same can be done for $\hat {\cal F}$ and ${\hat {\cal D}}$ starting from their root  $\frac 1 2$  with the sets of descendants  $F^{-1}(x)$ and $D^{-1}(x)$, respectively. 
Furthermore, we have the following commutative diagram
$$
\begin{CD}
J  @>\phi >>  I @> ? >> I \\
@VGVV  @VVFV @VVDV \\
J  @>\phi >>  I @> ? >> I \\
\end{CD}
$$
\end{theorem}
\vsni
\noindent
{\sl Proof.} The first assertion follows from Lemma \ref{desce},  eq. (\ref {inve}) and Lemma \ref{genera}. The proof of the conjugation between $G$ and $F$ is immediate. That for $F$ and $D$ can be obtained reasoning along the same lines as in the proof of Theorem \ref{conju}, starting from the observation that $D$ acts a the shift on binary expansions whereas the action of $F$ is the Farey shift $[0;a_1,a_2, \dots] \mapsto  [0;a_1-1,a_2, \dots]$ on the interval $[0, {1\over 2}]$ and  $[0;1,a_2, \dots] \mapsto  1-[0;a_2, \dots]$ on $({1\over 2} ,1]$. Then use Lemma \ref{reci}. We leave the details to the interested reader. $\qed$
\noindent
\begin{remark}{\rm 
Conversely, using the maps $G$, $F$ and $D$ one can retrace
the path from a leaf $x$ in any of the trees $\cal T$, $\cal F$ or $\cal D$ back to the root. For instance, for $x\in {\cal T}$ let $X=\prod_{i=1}^k M_i$ be the element which uniquely represents $x$ in $SL(2,\Z)$ with $k={\rm depth} (x)$, according to Proposition \ref{codaggio}. One then sees that the
following rule is in force: if ${G}^{(i-1)}(x)<1$ then $M_i=L$, ${G}^{(i-1)}(x)>1$ then $M_i=R$, for
$i=1,\dots, k$ with $k= {\rm depth} (x)$ so that ${G}^k(x)=1$. }
\end{remark}
\begin{remark}\label{rema} {\rm The map $D$ preserves the Lebesgue measure $dx$ on $I$, whereas the map $F$ preserves the a.c. infinite measure $\mu(dx)=dx/x(1-x)=\left( {\frac d  {dx}} \log \phi^{-1}(x)\right)dx$ on $I$, as one easily checks. This entails that $G$ preserves the (infinite) measure $\nu (dx) = \mu \circ \phi ( dx) = dx/x$ on $J$.

\noindent
Note that the entropy of $(I,F,d\mu)$ is zero (as well as that of $(J,G,d\nu)$). 
On the other hand, from the above  theorem it follows that 
also the measure $d?$ is invariant under $F$ (as well as $d\rho$ for $G$) and the
entropy of $(I,F,d?)$ is $\log 2$. Therefore $d?$ is the {\sl measure of maximal entropy} for $(I,F)$ (as well as $d\rho$ for $(J,G)$).}
\end{remark}
\noindent
To  the map $G$ we associate a {\sl generalised transfer operator} $L_q$ acting on $f:J\to \C$ as
\be\label{tranG}
(L_qf)(x)=\sum_{y\in G^{-1}(x)} \frac {f(y)} {|G'(y)|^q}
\ee
or else
\be
(L_qf)(x)=\frac 1 {(1+x)^{2q}}\, f\left(\frac x {1+x}\right) +  f(x+1)
\ee
where $q$ is a real or complex parameter. We point  out that a continuous fixed function for $L_q$  satisfies the functional equation
\be
f(x) =f(x+1)+ \frac 1 {(1+x)^{2q}}\, f\left(\frac x {1+x}\right) 
\ee
which is called {\sl Lewis-Zagier three-term functional equation} and is related to the spectral theory of the hyperbolic laplacian on the modular surface (see \cite{LeZa} and references therein). 

\noindent
In the same way,  the operators associated to $D$ and $F$ act on $f:I\to \C$ as
\be \label{dyatra}
f(x)\; \mapsto \;  \frac 1 {2^q} f\left(\frac x 2\right) +\frac 1 {2^q}  f\left(\frac x 2 + \frac 1 2 \right)
\ee
and
\be\label{fartra}
f(x)\; \mapsto \; \frac 1 {(1+x)^{2q}}\, f\left(\frac x {1+x}\right) + \frac 1 {(2-x)^{2q}} \, f\left(\frac 1 {2-x}\right)
\ee
respectively. For the spectral theory of an operator closely related to (\ref{fartra}) see \cite{I} and \cite{BGI}.

\noindent

\subsection{Harmonic functions and martingales}

\noindent
Let
$\Phi_s$, $s\in\{0,1\}$, be the inverse branches of $G$, i.e. \be\label{inve} \Phi_0(x)=\frac x {1+x}\quad ,\quad \Phi_1(x)=x+1
\ee
They satisfy:
\be \label{sym0}
\Phi_s(1/x)=\frac 1{\Phi_{1-s}(x)}\quad , \quad s\in \{0,1\}
\ee

\noindent
Let moreover $p(s,\cdot )$, $s\in \{0,1\}$, be a pair of positive Borel functions such that $p(0,x)+p(1,x)=1$, 
$\forall x\in J$.
We now want to study the Markov chain with state space $J$ where at each step, starting from a state $x \in J$, two transitions are possible towards the states $\Phi_0(x)$ and  $\Phi_1(x)$, with probabilities  $p(0,x)$ and $p(1,x)$ respectively.  Note that for $x=\frac 1 1$  and $p(i,x)=\frac 1 2$, $i=0,1$, this Markov chain  reduces to the random walk on ${\hat {\cal T}}$ discussed in Theorem \ref{rw}. 

\noindent
We now briefly adapt to our context some basic facts about canonical  Markov chains associated to Markov transfer operators (see \cite{CoRa}; also \cite{CoRa1} for an application to the dyadic transfer operator (\ref{dyatra})). 
 Let $P: L^{\infty}(J) \to L^{\infty}(J)$ be the Markov operator acting as
\be
(P f)(x)=p(0,x) \, f\left(\Phi_0(x)\right) + p(1,x) \, f\left(\Phi_1(x)\right)
\ee
A measurable function $h:J\to \C$ satisfying $Ph=h$ is called $P$-{\sl harmonic}.
In the sequel we shall make the further assumption that the transition probabilities 
satisfy:
\be \label{sym1}
p(s, 1/x) = p(1-s, x) \quad , \quad s\in\{0,1\}, \quad \forall x\in J
\ee
The symmetries (\ref{sym0}) and  (\ref{sym1}) yield at once the following
\begin{lemma} \label{commuta} The  averaging operator $A: L^{\infty}(J) \to L^{\infty}(J)$ acting as
\be
(Af)(x)=\frac {f(x)+f(1/x)} 2
\ee
commutes with $P$. In particular, if $h:J\to \C$ is a bounded $P$-harmonic function then $Ah$ has the same property.
\end{lemma}
\noindent
A positive measure $\nu$ is called $P$-{\sl invariant} if $\nu P = \nu$, i.e. $\int_J Pf \, d\nu = \int_J f \, d\nu$ for all measurable $f:J\to \C$. In turn, one readily realizes that this condition  is equivalent to
\be
\frac {d \nu \circ \Phi_s} {d \nu} = p(s,\cdot) , \quad s\in \{0,1\}
\ee

\noindent
Now, setting $\Omega :=\{0,1\}^{\N}$, a $n$-dimensional {\sl cylinder} of $\Omega$ is a subset of the type $C(i_1, \dots, i_n)=\{ \omega \in \Omega \, :| \, \omega_1=i_1, \dots, \omega_n = i_n\}$. The cylinder sets generate the topology of $\Omega$ and its Borel $\sigma$-algebra $\cal F$. 

\noindent
Given $x\in J$ let $U(x)$ be the closure of the set of all possible paths starting at $x$, i.e.
\be
U(x)= {\overline {\cup_{\omega \in \Omega }
\{\Phi_{\omega_{n}}\circ \cdots \circ \Phi_{\omega_{1}}(x)\, , \, n\geq 1\} }}
\ee
This is clearly a compact invariant set, in the sense that if $y\in U(x)$ then $\Phi_i(y)\in U(x)$, $i\in \{0,1\}$.
More generally, a compact subset $V$ of $J$ is called {\sl invariant} if for all $x\in V$ and all $i\in \{0,1\}$ such that $p(i,x)>0$ we have $\Phi_i(x)\in V$.

\noindent
A first basic fact (see \cite{CoRa}, Sec. 3.4; or else \cite{Jo}, Chap. 2.4) is that  for each $x \in J$  there is a unique probability measure $\,{\Ps}_x$ on $\Omega$ such that
\be 
{\Ps}_x(C(i_1, \dots, i_n))=\prod_{k=1}^n p(i_k,\Phi_{i_{k-1}}\circ \cdots \circ \Phi_{i_{1}}(x))
\ee
\vsni
\noindent
The symmetries (\ref{sym0}) and (\ref{sym1}) entail the following
\begin{lemma} \label{symme} For each $x \in J$ we have
\be
{\Ps}_x(C(i_1, \dots, i_n))={\Ps}_{1/x}(C(1-i_1, \dots,1- i_n))
\ee
\end{lemma}
\vsni
\noindent
For $\omega \in \Omega$ let $X_k(\omega)=\omega_k$ be the $k$-th coordinate function on $\Omega$ and ${\cal Z}_n$ the subalgebra of $C(\Omega)$ generated by the first $n$ coordinates $\{X_k, 1\leq k\leq n\}$. The ${\cal Z}_n$ form a filtration in that ${\cal Z}_n \subset {\cal Z}_{n+1}$. 
Note that if $X\in {\cal Z}_n$ then 
$$
{\Es}_x [X] = \sum_{(\omega_1, \cdots,  \omega_n)\in \{0,1\}^{n} }
\prod_{k=1}^n 
p({\omega_k},\Phi_{\omega_{k-1}}\circ \cdots \circ \Phi_{\omega_{1}}(x))\, 
X(\omega_1, \cdots , \omega_n)
$$
In particular, if there is  $h:J\to \C$ s.t. $$X(\omega_1, \cdots , \omega_n)= h(\Phi_{\omega_n}\circ \cdots \circ \Phi_{\omega_{1}}(x))$$ then
\be\label{ide}
{\Es}_x [X] = (P^n h)(x)
\ee

\vsni
\noindent
Now, having fixed $x\in J$, define
\be\label{mc}
W_0(x,\omega):=x\quad \hbox{and} \quad W_n(x,\omega):= \Phi_{X_n(\omega)}\circ \cdots \circ \Phi_{X_1(\omega)}(x),\quad n\geq 1
\ee
The process $\{ W_n(x,\,\cdot \,), \, n\geq 0\}$ defined on $(\Omega, {\cal F}, {\Ps}_x)$ is a Markov chain on $J$ with initial state $x$ and for any measurable function $f:J\to \C$ we have
\be
\lim_{n\to \infty} (P^nf)(x)=\lim_{n\to \infty} {\Es}_x [f(W_n(x,\, \cdot \,)]
\ee

\vsni
\noindent
Moreover, if $h:J\to \C$ is a measurable bounded $P$-harmonic function then we have
\begin{eqnarray}
{\Es}_x [ h(W_{n+1}(x,\,\cdot \,)) \, | \,   {\cal Z}_{n}  ]  &=& \sum_{\omega_{n+1} \in \{0,1\}} p(\omega_{n+1},x) 
h(\Phi_{\omega_{n+1}}\circ \cdots \circ \Phi_{\omega_{1}}(x)) \nonumber \\
&=& (Ph)(\Phi_{\omega_{n}}\circ \cdots \circ \Phi_{\omega_{1}}(x)) \nonumber \\
&=& h(\Phi_{\omega_{n}}\circ \cdots \circ \Phi_{\omega_{1}}(x)) \nonumber \\
&=&  h(W_{n}(x,\,\cdot \,)) \nonumber
\end{eqnarray} 
In other words the sequence of random variables $\{ h(W_n(x,\,\cdot \,)), \, n\geq 0\}$ on $(\Omega, {\cal F}, {\Ps}_x)$  is a bounded martingale (relative to the filtration $\{{\cal Z}_n, n\geq 1\}$) and therefore it converges pointwise ${\Ps}_x$-a.e. The limit random variable $H(x,\,\cdot \,)=\lim_{n\to \infty} h(W_n(x,\,\cdot \,))$ satisfies
\be\label{cocycle}
H(x,\omega)=H(\Phi_{\omega_1}(x), \sigma \omega)
\ee
where $\sigma : \Omega \to \Omega$ is the left shift acting as $(\sigma \omega)_i=\omega_{i+1}$.  A bounded measurable function $H:J\times \Omega \to \C$ satisfying (\ref{cocycle}) is said to be a {\sl cocycle}. Conversely, by (\ref{ide}) $h$ may be recovered from the cocycle $H$ as
\be \label{boundary}
h(x) = {\Es}_x [H(x,\,\cdot \,)]
\ee
\begin{remark} Note that $P: L^{\infty}(J) \to L^{\infty}(J)$ has norm one. Therefore if $h\in L^{\infty}(J)$ is an eigenfuction of $P$ corresponding to a real and positive eigenvalue then the sequence $\{ h(W_n(x,\,\cdot \,)), \, n\geq 0\}$ on $(\Omega, {\cal F}, {\Ps}_x)$  is a supermartingale, which again converges ${\Ps}_x$-a.e. to a limit cocycle $H$.
\end{remark}
\begin{remark} As pointed out in \cite{Jo}, p.50, eq. (\ref{boundary}) can be thought of as an analogue of the classical result about the existence of boundary functions for bounded harmonic functions via Poisson integral.
\end{remark}

\noindent
We now discuss two specific Markov chains of the above type, denoted $MC^{(0)}$ and  $MC^{(1)}$, corresponding to the choices $q=0$ and $q=1$ in (\ref{tranG}). 
\subsubsection{The Markov chain $MC^{(0)}$}

\noindent
Setting $q=0$ in  (\ref{tranG}) we have ${\frac 1 2}L_0 1 = 1$. One can then consider the Markov (i.e. normalised) operator $P^{(0)}$ acting as 
$P^{(0)} f={\frac 1 2}L_0f$. More esplicitly, 
\be\label{P^{(0)}}
(P^{(0)} f)(x)=\frac 1 2 \, f\left(\frac x {1+x}\right) + \frac 1 2\, f\left(x+1\right)
\ee
Note that if $h$ is $P^{(0)}$-harmonic then, iterating (\ref{P^{(0)}}) we get
$$
h(x)=\sum_{k=0}^{N-1} \frac 1 {2^{k+1}} \, h\left( \frac {x+k}{x+k+1}\right) + \frac 1 {2^{N}} \, h(x+N)
$$
We therefore have the
\begin{lemma} 
A bounded function $h:J \to \C$ is $P^{(0)}$-harmonic if and only if 
$$
h(x)=\sum_{k=0}^{\infty} \frac 1 {2^{k+1}}\, 
h\left( \frac {x+k}{x+k+1}\right)\quad , \quad x\in [0,\infty)
$$
\end{lemma}

\noindent
Moreover,
\vsni
\noindent
\begin{lemma} Let $\rho$ be as in (\ref{defi}).
The probability measure $d\rho$ on $J$  is $P^{(0)}$-invariant.
\end{lemma}
\noindent
{\sl Proof.} From the fact that the function $\rho$ is the distribution function of the (permuted) Stern-Brocot fractions (cf. Theorem \ref{distri} ) and Lemma \ref{desce}  one readily obtains that $\rho$ satisfies the functional equation 
\be
2\rho (x) =  \left\{
\begin{array}{cl}
\displaystyle \rho\left(\frac x {1-x}\right)  & {\rm if}\ 0< x< {1} \\[0.5cm]
\rho(x-1)+1 & {\rm if}\  x\geq 1
\end{array} \right.  
\ee
The claim now follows straightforwardly. $\qed$
\vsni
\noindent
Setting $p(0,x)=p(1,x)=1/2$ we have that there are no compact invariant sets and according to (\cite{CoRa1}, Sec. IV) ) $h\equiv 1$ is the only bounded continuous $P^{(0)}$-harmonic function. 
Moreover, the 
unique probability measure $\,{\Ps}_x^{(0)}$  on $\Omega$ such that
\be 
{\Ps}_x^{(0)}(C(i_1, \dots, i_n))=2^{-n}
\ee
 is atomless for each $x\in J$. The Markov chain $MC^{(0)}$ is then defined as in (\ref{mc}) on the probability space $(\Omega, {\cal F}, {\Ps}_x^{(0)})$. We summarize the above in the following
 \begin{theorem}
 For $f\in L^1(\R^+, d\rho )$ we have
 \be
\lim_{n\to \infty} ({P^{(0)}}^nf)(x)=\lim_{n\to \infty} {\Es}^{(0)}_x [f(W_n(x,\, \cdot \,)]=\int_0^\infty f \, d\rho
\ee
\end{theorem}

\noindent
Taking $f=1_{(a,b)}$, $(a,b)\subset \R^+$,  this is to be compared  with Theorem \ref{rw}. 
\subsubsection{The Markov chain $MC^{(1)}$}

Setting $q=1$ in  (\ref{tranG}) we have $L_1 g=g$ where $g(x)=1/x$ is the $G$-invariant density. 
We then consider the Markov operator $P^{(1)}$ acting as 
$P^{(1)}f=g^{-1}L_1(f \cdot g)$, or 
\be\label{fixpt}
(P^{(1)} f)(x)=\frac 1 {x+1} \, f\left(\frac x {1+x}\right) + \frac x {1+x}\, f\left(x+1\right)
\ee
If $h$ is  $P^{(1)}$-harmonic then
$$
h(x)=\sum_{k=0}^{N-1} \frac x {(x+k)(x+k+1)} \, h\left( \frac {x+k}{x+k+1}\right) + \frac x {x+N} \, h(x+N)
$$
\begin{lemma} 
A bounded function $h:J \to \C$ is $P^{(1)}$-harmonic if and only if 
$$
h(x)=\sum_{k=0}^{\infty} \frac x {(x+k)(x+k+1)}\, 
h\left( \frac {x+k}{x+k+1}\right)\quad , \quad x\in [0,\infty)
$$
\end{lemma}
\noindent
Furthermore, the validity of $L_1g=g$ is equivalent to the fact that the infinite measure $\nu(dx)=dx/x$ on $J$ is $P^{(1)}$-invariant.

\noindent
Set moreover
\be\label{probba1}
p(0,\infty)=p(1,0)=0\quad , \quad p(0,0)=p(1,\infty)=1
\ee
and
\be\label{probba11} p(0,x)=\frac 1 {x+1}\quad ,\quad p(1,x)=\frac x {x+1}\quad , \quad x\in (0,\infty) \ee
They plainly satisfy the symmetry (\ref{sym1}).
\noindent
Moreover, from (\ref{probba1}) it follows that the singletons $\{0\}$ and $\{\infty\}$ are two {\sl disjoint compact invariant sets} and from (\ref{probba11}) one sees that they are the only invariant sets of this type.
\vsni
\noindent
The Markov chain $MC^{(1)}$ is now defined as in (\ref{mc}) on the probability space $(\Omega, {\cal F}, {\Ps}_x^{(1)})$, where  ${\Ps}_x^{(1)}$ is the transition measure  on $\Omega$ arising from the probabilities (\ref{probba1}) and (\ref{probba11}). It satisfies the following
\begin{lemma} \label{atomic}
For $x\in (0,\infty)$ the measures $\,{\Ps}_x^{(1)}$ have no atoms.  On the other hand, both $\,{\Ps}^{(1)}_0$ and $\,{\Ps}^{(1)}_\infty$ are purely atomic with  $\,{\Ps}^{(1)}_0= \delta_{0^\infty}$ and   $\,{\Ps}^{(1)}_\infty= \delta_{1^\infty}$.
\end{lemma}
\noindent
{\sl Proof.} From (\ref{inve}) and (\ref{probba11}) it follows that the path of length $n$ starting at $x\in [1,\infty)$ and having largest probability is that corresponding to the word $\omega=1\cdots 1$. If instead $0<x<1$ it corresponds to $\omega=0\cdots 0$. On the other hand we have
$$
{\Ps}^{(1)}_x(C(1, \dots, 1)) ={\Ps}^{(1)}_{1/x}(C(0, \dots, 0)) = \prod_{k=0}^{n-1} \frac {x+k}{x+k+1} = \frac x {x+n} \to 0
$$
as $n\to \infty$, proving the first assertion. The last is straigthforward. 
$\qed$ 
\vsni
\noindent
Now, a path starting somewhere in $J$ and converging to $0$ corresponds to a sequence of the form
$(\omega_1,\dots, \omega_n, 0, 0, 0,  \dots)$ for some $n\geq 1$. The symmetric sequence $(1-\omega_1,\dots, 1-\omega_n, 1, 1, 1,  \dots)$ yields a corresponding path which converges to $\infty$. 
By Lemma \ref{symme}, if we let the first path start at $x$ and the second one at $1/x$, all finite equal
portions of them have the same probability. We can thus concentrate on the paths starting at $x$ and converging to $0$. In turn, these can be put in a one-to-one correspondence with 
$\Q_2$ via the mapping 
$$
\Q_2 \ni a =\sum_{i=1}^n \omega_i 2^{-i} \mapsto \omega(a)=(\omega_1,\dots, \omega_n, 0, 0, 0,  \dots)
$$
Another copy of $\Q_2$ is obtained via the mapping
$$
1-\Q_2 \ni 1-a =\sum_{i=1}^n (1-\omega_i) 2^{-i} +2^{-n}\mapsto \omega(1-a)=(1-\omega_1,\dots, 1-\omega_n, 1, 1, 1,  \dots)
$$
With the identification $a \leftrightarrow \omega (a)$ we set
\be 
{\Ps}^{(1)}_x(\Q_2) = \sum_{a\in \Q_2} {\Ps}^{(1)}_x( \omega (a))
\ee
so that Lemma \ref{atomic} can be rephrased in the form
\be\label{cudue}
{\Ps}^{(1)}_x(\Q_2) +{\Ps}^{(1)}_x(1-\Q_2) = \delta_x^0 +\delta_x^\infty
\ee

\vsni
\noindent

\noindent
Finally, putting together the above and  (\cite{CoRa1}, Sec. IV)  we get the following
\begin{theorem} 
The space of bounded harmonic functions has dimension two. 
A basis for it is given by the functions 
$$
h_0(x)=1 \quad \hbox{and}\quad  h_1(x) = w(x)+w(1/x)
$$
where 
$$
 w(x)={\Ps}^{(1)}_x \left[ \lim_{n\to \infty} W_n(x,\,\cdot \,) =0\right]=  {\Ps}^{(1)}_{1/x} \left[ \lim_{n\to \infty} W_n\left( 1/x ,\,\cdot \,\right) =\infty\right]
$$
\end{theorem}
\noindent
\begin{remark}
Having fixed $x,y\in J$ let $B(x,y)$ be the tail event
$$
B(x,y) := \{ \omega \in \Omega \, : \,  \lim_{n\to \infty} W_n(x, \omega) =y\} 
$$
According to (\ref{boundary}),  the cocycle $H_1$ associated to  $h_1$ is
$$
H_1(x,\omega) = 1_{B(x,0)\cup B(x,\infty )}(\omega)
$$
Note that, eq. (\ref{cudue}) can be further rephrased as
\be
{\Ps}^{(1)}_x \left[ B(x,0)\right]+ {\Ps}^{(1)}_x \left[B(x,\infty )\right]=\delta_x^0 +\delta_x^\infty
\ee
\end{remark}

\end{document}